\documentclass{amsart}
\title{Parking functions and Haglund--Loehr data}
\author[Yu.Burman]{Yurii M. Burman}
\address{Independent University of Moscow, 121002, 11, B.Vlassievsky per., 
Moscow, Russia}
\email{burman@mccme.ru}
\thanks{Author was supported by the INTAS grant YSF-2001/1-81 and by the 
RFBR grants Nos 01-01-00660 and N.Sh.1972.2003.1.}
\date{}

\makeatletter

\RequirePackage{amssymb}

\newcommand{\theoremName}{Theorem}
\newcommand{\lemmaName}{Lemma}
\newcommand{\corollaryName}{Corollary}
\newcommand{\statementName}{Statement}
\newcommand{\remarkName}{Remark}
\newcommand{\exampleName}{Example}
\newcommand{\definitionName}{Definition}
\newcommand{\problemName}{Problem}
\newcommand{\proofName}{Proof}
\renewcommand{\proofname}{\proofName}
\newcommand{\answerName}{Answer}
\newcommand{\hintName}{Hint}

\theoremstyle{plain}

\newtheorem {Statement}{\statementName}

\theoremstyle{remark}

\theoremstyle{definition}

\let\@newpf\proof \let\proof\relax 
\def \namepf[#1] {\@newpf[\proofname\ #1]}
\newenvironment{proof}{\@ifnextchar[{\namepf}{\@newpf[\proofname]}}{\qed\endtrivlist}

\def \Complex {{\mathbb C}}

\def \lnorm#1\rnorm {\vphantom{#1}\left\|\smash{#1}\right\|}
\def \lmod#1\rmod {\vphantom{#1}\left|\smash{#1}\right|}

\newcommand \bydef {\stackrel{\mbox{\scriptsize def}}{=}}

\@ifundefined{name}{\newcommand \name[1] {\mathop{\rm #1}\nolimits}}%
{\relax}

\renewcommand \phi {\varphi}
\renewcommand \rho {\varrho}

\makeatother

\begin{document}

 \begin{abstract}
A {\em parking function} is a sequence of $n$ nonnegative integers 
majorated by a permutation of the set $\{0, \dots, n-1\}$. We provide a way 
to encode parking functions by data suggested by J.\,Haglund and N.\,Loehr 
in \cite{HL}. This coding is compared with another one proposed earlier by 
M.\,Shapiro and the author.
 \end{abstract}

\maketitle

\section{The correspondence}

Denote $\Sigma_n$ the set of all permutations of $\{0,1,\dots,n-1\}$; 
elements of $\Sigma_n$ will be simply referred to as ``permutations''. A 
number $t \in \{0,1,\dots,n-1\}$ is called a {\em descent} of $\sigma$ if 
$\sigma(t-1) > \sigma(t)$. The sum 
 \begin{equation*}
\name{maj}(\sigma) = \sum_{\text{$t$ is a descent of $\sigma$}} (n-t)
 \end{equation*}
is called the major index of $\sigma$. 

In \cite{HL} J.\,Haglund and N.\,Loehr introduced the following 
characteristics of a permutation $\sigma$. Extend the sequence $\sigma(0), 
\dots, \sigma(n-1)$ by the term $\sigma(-1) \bydef +\infty$. Then for every 
$i = 0, \dots, n-1$ define:

\vskip\parskip\hangindent\parindent {\em
$u_i(\sigma)$ is the smallest $j,\, 0 \le j \le i$, such that the segment 
$\sigma(j-1), \sigma(j), \dots, \sigma(i)$ either contains no descents or 
contains exactly one descent and $\sigma(j-1) > \sigma(i)$. 
}\vskip\parskip

\noindent The sequence $u_0(\sigma), \dots, u_{n-1}(\sigma)$ will be 
denoted $u(\sigma)$. It is easy to see that $u_0(\sigma) \le u_1(\sigma) 
\le \dots \le u_{n-1}(\sigma)$ and $0 \le u_i(\sigma) \le i$ for all $i = 
0, \dots, n-1$.

A pair $(\sigma,k)$ will be called an HL-pair if $\sigma$ is a permutation, 
and $k = (k_0, \dots, k_{n-1})$ is a sequence of integers such that 
$u_i(\sigma) \le k_i \le i$ for all $i = 0, \dots, n-1$. A sequence of 
nonnegative integers $p_0, \dots, p_{n-1}$ is called a {\em parking 
function} if it is majorated by a permutation $\tau \in \Sigma_n$: $0 
\le p_i \le \tau(i)$ for all $i$.

In \cite{HL} the authors consider the following generating function:
 \begin{equation*}
R_n(q,t) = \sum_{\text{$(\sigma,k)$ is an HL-pair}} q^{\name{maj}(\sigma)} 
t^{n(n-1)/2 - \sum_i k_i}
 \end{equation*}
and express it as a certain sum over the set of parking functions. They 
also conjecture that the coefficient at $q^at^b$ in $R_n(q,t)$ equals to 
the dimension of the bihomogenous component $H_{a,b}$ of degree $(a,b)$ in 
the module of diagonal coinvariants. This module is defined as a quotient 
of the ring $\Complex[x_1, \dots, x_n, y_1, \dots, y_n]$ by the ideal $J$ 
generated by all the polynomials $P$ without a constant term invariant 
under the diagonal action of $\Sigma_n$: $\sigma P = P$ where $\sigma 
P(x_1, \dots, x_n, y_1, \dots, y_n) \bydef P(x_{\sigma(1)}, \dots, 
x_{\sigma(n)}, y_{\sigma(1)}, \dots, y_{\sigma(n)})$. A degree of a 
diagonal coinvariant is a pair $(a,b)$ where $a$ is its total degree with 
respect to all the $x_i$ and $b$, a degree with respect to all the $y_i$. 
See \cite{HInit} for details.

Here we present a direct one-to-one correspondence between the sets of 
HL-pairs and parking functions. Namely, for an HL-pair $(\sigma,k)$ denote 
$p(\sigma,k)$ the sequence $p_0, \dots, p_{n-1}$ such that $p_{\sigma(i)} = 
k_i$ for all $i \in I_n$. 

 \begin{Statement}
The mapping $(\sigma,k) \mapsto p(\sigma,k)$ is a one-to-one correspondence 
between the sets of HL-pairs and parking functions.
 \end{Statement}

 \begin{proof}
Since $k_i \le i$, one has $p(\sigma,k)_i \le \sigma^{-1}(i)$, and 
therefore $p(\sigma,k)$ is a parking function. We are going to prove that 
for every parking function $q = (q_0, \dots, q_{n-1})$ there exists exactly 
one HL-pair $(\sigma,k)$ such that $q = p(\sigma,k)$.

Consider {\bfseries Case I} when the parking function $q = (q_0, \dots, 
q_{n-1})$ is nondecreasing: $q_0 \le q_1 \le \dots \le q_{n-1}$. Take
$\sigma$ to be the identity permutation. Then $u_i(\sigma) = 0$ for all 
$i$, and $(e,k)$ is an HL-pair for every sequence $k$ such that $0 \le 
k_i \le i$ for all $i$. One has now $q = p(e,q)$, and the existence of 
$(\sigma,k)$ for such $q$ is proved. 

To prove uniqueness, let $q = p(\sigma,k)$. For a nondecreasing parking 
function $q$ one has $p_i \le i \Longleftrightarrow k_i \le \sigma(i)$ for 
all $i$, and therefore $u_i(\sigma) \le \sigma(i)$. Take $i = 
\sigma^{-1}(0)$, so that one has $u_i(\sigma) = 0$. Since $\sigma(-1) = 
+\infty > \sigma(i) = 0$, it is possible only if $\sigma(0), \sigma(1), 
\dots, \sigma(i) = 0$ contains no descents, so that $i = 0$. Now make an 
induction step: suppose $\sigma(j) = j$ for all $j = 0, \dots, t-1$, and 
let $i = \sigma^{-1}(t)$. Then $a \bydef u_i(\sigma) \le t$, and similar to 
the previous reasoning one concludes that $\sigma(a) = a, \dots, 
\sigma(t-1) = t-1. \dots,  \sigma(i) = t$ must be a sequence without 
descents, hence $i = t$. So we proved by induction that $\sigma$ is an 
identity permutation. Now $q = p(\sigma,k)$ implies that $k = q$, and the 
uniqueness in Case I is proved.

{\bfseries Case II}: general. For a parking function $q = (q_0, \dots, 
q_{n-1})$ denote $\ell(q) = \sum_{i = 0}^{n-1} (n-i) q_i$. We are going to 
prove the existence and uniqueness of $(\sigma,k)$ by induction on 
$\ell(q)$. Case I constitutes the induction base.

First, prove the existence. If $q$ is nondecreasing, it is Case I, so let 
$t$ be such that $q_{t-1} > q_t$. Construct a parking function $q'$ such 
that $q'_{t-1} = q_t$, $q'_t = q_{t-1}$, and $q'_i = q_i$ for all the other 
$i$. One has $\ell(q') < \ell(q)$, and by the induction hypothesis, there 
exists an HL-pair $(\sigma,k)$ such that $p(\sigma,k) = q'$. Let $i = 
\sigma^{-1}(t-1)$ and $j = \sigma^{-1}(t)$; consider now the following 
cases.

\subsection{}\label{It:Distant} $j \ge i+2$. Define a permutation $\sigma'$ 
such that $\sigma'(i) = t$, $\sigma'(j) = t-1$ and $\sigma'(a) = \sigma(a)$ 
for all the other $a$. Apparently, $u(\sigma') = u(\sigma)$ and therefore 
$(\sigma',k)$ is an HL-pair. One has $p(\sigma',k) = q$, so the existence 
is proved.

\subsection{} $j = i+1$. In this case $u_j(\sigma) = u_i(\sigma)$, and the 
case splits in two:

\subsubsection{}\label{It:Adjacent} $p_t = k_j \le i$. Define a sequence 
$k'$ as $k'_i = k_{i+1}$, $k'_{i+1} = k_i$, and $k'_a = k_a$ for all the 
other $a$. Apparently, $(\sigma,k')$ is an HL-pair and $p(\sigma,k') = q$.

\subsubsection{} $p_t = k_j = i+1$. Define a permutation $\sigma'$ as in 
Case \ref{It:Distant}. Now $u_a(\sigma') = u_a(\sigma)$ for all $a \ne 
i+1$, and $u_{i+1}(\sigma') = i+1$, so $(\sigma',k)$ is an HL-pair. 
Apparently, $p(\sigma',k) = q$. 

\subsection{} $j < i$. Since the values $\sigma(j) = t$ and $\sigma(i) = 
t-1$ are adjacent, one has $u_i(\sigma) \ge j+1$, and therefore $p_{t-1} = 
k_i \ge u_i(\sigma) \ge j+1 > j \ge k_j = p_t$, contrary to the assumption. 
So, this case is impossible.

Thus, existence of the HL-pair $(\sigma,k)$ such that $p(\sigma,k) = q$ is 
proved for every $q$.

Prove now the uniqueness. Take a parking function $q = p(\sigma,k)$ and 
choose $t$ such that $q_{t-1} > q_t$. Let $i = \sigma^{-1}(t-1)$ and $j = 
\sigma^{-1}(t)$. Define a parking function $q'$ as before ($q_t$ and 
$q_{t-1}$ are exchanged). Now we will fix a specific HL-pair $(\sigma',k')$ 
such that $p(\sigma',k') = q'$ (since the uniqueness is not proved yet, 
there could be several possibilities here). To do this consider the 
following cases:

\setcounter{subsection}{0}
\subsection{} $j \ge i+2$. Then $u_j(\sigma) \ge i+2$, and therefore $q_t = 
k_j \ge u_j(\sigma) \ge i+2 > i \ge k_i = q_{t-1}$, contrary to the 
assumption. So, this case is impossible.

\subsection{}\label{It:AdjPlus} $j = i+1$. Here $k_{i+1} = q_t < q_{t-1} = 
k_i \le i$ and $u_{i+1}(\sigma) = u_i(\sigma)$. Define $k'$ as in Case 
\ref{It:Adjacent} of the existence proof, and note that $(\sigma, k')$ is 
an HL-pair with $p(\sigma, k') = q'$.

\subsection{}\label{It:AdjMinus} $j = i-1$. Here $u_i(\sigma) = i$ and 
therefore $p_{t-1} = k_i = i$. Define $\sigma'$ as in Case \ref{It:Distant} 
of the existence proof. Apparently, $(\sigma',k)$ is an HL-pair and 
$p(\sigma',k) = q'$. 

\subsection{}\label{It:Minus} $j \le i-2$. Here, again, define $\sigma'$ as 
before, and note that $u(\sigma') = u(\sigma)$. This means that 
$(\sigma',k)$ is an HL-pair, and $p(\sigma',k) = q'$. 

Let now $q = p(\sigma^{(1)},k^{(1)}) = p(\sigma^{(2)},k^{(2)})$. If 
$\sigma^{(1)} = \sigma^{(2)}$ then $k^{(1)} = k^{(2)}$, so suppose that 
$\sigma^{(1)} \ne \sigma^{(2)}$. Apply the procedure just described to 
$(\sigma^{(1)},k^{(1)})$ and $(\sigma^{(2)},k^{(2)})$; it gives two 
representations: $q' = p((\sigma^{(1)})',(k^{(1)})') = 
p((\sigma^{(2)})',(k^{(2)})')$. Since $\ell(q') < \ell(q)$, these two 
representations must be the same by the induction hypothesis. This is 
possible only if the pair $(\sigma^{(1)}, k^{(1)})$ belongs to Case 
\ref{It:AdjPlus} above (i.e.\ $(\sigma^{(1)})^{-1}(t) = 
(\sigma^{(1)})^{-1}(t-1) + 1$ and $\sigma^{(1)} = (\sigma^{(1)})'$ but 
$k^{(1)} \ne (k^{(1)})'$), while the pair $(\sigma^{(2)}, k^{(2)})$ belongs 
to Case \ref{It:AdjMinus} or Case \ref{It:Minus} (so that $k^{(2)} = 
(k^{(2)})'$ but $\sigma^{(2)} \ne (\sigma^{(2)})'$). Hence, 
$(\sigma^{(2)})' = \sigma^{(1)}$ and $k^{(2)} = (k^{(1)})'$. If 
$\sigma^{(1)}(i_1) = \sigma^{(2)}(i_2) = t-1$ and $\sigma^{(1)}(j_1) = 
\sigma^{(2)}(j_2) = t$ then $i_1 = j_2 = i_2 - 1 = j_1 - 1$. This implies 
$u_{i_2}(\sigma^{(2)}) = i_2$ and therefore $k^{(2)}_{i_2} = i_2$. But then 
$k^{(1)}_{i_1} = k^{(2)}_{i_2} = i_2 = i_1 + 1$ which is impossible. The 
uniqueness is proved.
 \end{proof}

\section{HL-pairs $(\sigma,k)$ with small $\name{maj}(\sigma)$}

Here we investigate a connection between two constructions of the parking 
functions --- the one described above and the construction of the paper 
\cite{BS}. The latter is as follows: call a pair of integer sequences 
$(k_0, \dots, k_{n-1})$ and $(l_0, \dots, l_{n-1})$ {\em admissible} if 
 \begin{enumerate}
\item $0 \le l_i \le k_i \le i$ for all $i = 0, \dots, n-1$.
\item If $i < j$ and $l_i > l_j$ then $k_j \ge i+1$.
 \end{enumerate}
Define a permutation $\sigma_{k,l}$ such that $k_i - l_i = \#\{j < i \mid 
\sigma_{k,l}(j) > \sigma_{k,l}(i)\}$ for every $i = 0, \dots, n-1$ (it is 
easy to see that $\sigma_{k,l}$ exists and is unique), and then define a 
parking function $q = f(k,l)$ like before: $q_{\sigma_{k,l}(i)} = k_i$ for 
all $i = 0, \dots, n-1$. It is proved in \cite{BS} that $f$ provides a 
one-to-one correspondence between the set of admissible pairs and the set 
of parking functions. There was also conjectured in \cite{BS} that there 
exists an integer-valued function $b = b(l)$ such that the number of 
admissible pairs $(k,l)$ with $\sum_i k_i = n(n-1)/2 - a$ and $b(l) = b$ 
equals the dimension of a bihomogenous component $H_{a,b}$ of degree 
$(a,b)$ in the module of diagonal coinvariants. The paper \cite{BS} 
provides a conjectural description of the sets $b^{-1}(t)$ where $0 \le t 
\le 3$.

Compare this conjecture with the conjecture of \cite{HL} mentioned above. 
Let $(\sigma,k)$ be an HL-pair such that $\name{maj}(\sigma) \le 3$. 

\setcounter{subsection}{-1}
\subsection{} Let $\name{maj}(\sigma) = 0$, that is, $\sigma$ is an 
identity permutation $e$. The pair $(e,k)$ is HL for every $k_0, \dots, 
k_{n-1}$ such that $0 \le k_i \le i$ for all $i = 0, \dots, n-1$. It is 
easy to see that the number of such pairs equals to the dimension of 
$H_{a,0}$ where $a = n(n-1)/2 - \sum_i k_i$ (so, the conjecture of 
\cite{HL} is true in this case). 

To every HL-pair $(e,k)$ we associate a pair $(k,0)$ (where $0$ denotes a 
sequence $l$ of $n$ zeros). This pair is admissible. So, the dimension of 
$H_{a,0}$ equals to the number of admissible pairs $(k,l)$ such that 
$\sum_i k_i = n(n-1)/2 - a$ and $l = 0$. This coincides with the conjecture 
of \cite{BS} saying that $b^{-1}(0) = \{0\}$.

\subsection{} Let $\name{maj}(\sigma) = 1$. This means that $\sigma$ has 
the form $(0,1, \dots, s-1, s+1, \dots, n-1, s)$ for some $s \le n-2$. One 
has $u_i(\sigma) = 0$ for $0 \le i \le n-2$, and $u_{n-1}(\sigma) = s+1$. 
Relate to an HL-pair $(\sigma,k)$ a pair $(k,l)$ where $l_i = 0$ for $0 \le 
i \le n-2$ and $l_{n-1} = s+1$. It is easy to see that $(k,l)$ is 
admissible and $f(k,l) = p(\sigma,k)$. Vice versa, if $(k,l)$ is an 
admissible pair with $l_0 = \dots = l_{n-2} = 0$ and $1 \le l_{n-1} \le 
n-1$ then $(\sigma,k)$ is an HL-pair (where $\sigma$ is described by the 
formula above, with $s = l_{n-1}-1$). Thus, according to the conjecture of 
\cite{HL}, the dimension of $H_{a,1}$ equals to the number of admissible 
pairs with $\sum_i k_i = n(n-1)/2-a$ and $l$ such that all its terms are 
zeros except the last one. This coincides with the conjecture of \cite{BS} 
saying that $b(l) = 1$ exactly for such $l$. In fact, both conjectures 
follow from the results of \cite{Alfano}.

\subsection{} Let $\name{maj}(\sigma) = 2$. In this case $\sigma$ has the 
form $(0, 1, \dots, s_1-1, s_1 + 1, \dots, s_2 - 1, s_2 + 1, \dots, n-1, 
s_1, s_2)$ where $s_1 < s_2$ and $s_1 \le n-3$. Therefore $u_i(\sigma) = 0$ 
for $0 \le i \le n-3$, and $u_{n-2}(\sigma) = s_1 + 1$, $u_{n-1}(\sigma) = 
s_2$. Relate to an HL-pair $(\sigma,k)$ a pair $(k,l)$ where $l_i = 0$ for 
$i \le n-2$, $l_{n-2} = s_1 + 1$ and $l_{n-1} = s_2$. It is easy to see 
that $(k,l)$ is admissible and $f(k,l) = p(\sigma,k)$. Vice versa, if 
$(k,l)$ is an admissible pair with $l$ such that $l_0 = \dots = l_{n-3} = 
0$ and $1 \le l_{n-2} \le l_{n-1}$, one can construct an HL-pair 
$(\sigma,k)$ by the formula above with $s_1 = l_{n-2}-1$ and $s_2 = 
l_{n-1}$. Thus, according to the conjecture of \cite{HL}, the dimension of 
$H_{a,2}$ equals to the number of admissible pairs with $\sum_i k_i = 
n(n-1)/2-a$ and $l$ just described. This coincides with the conjecture of 
\cite{BS} saying that $b(l) = 2$ exactly for such $l$.

\subsection{} Let $\name{maj}(\sigma) = 3$. In this case $\sigma$ has one 
of the following forms:

\subsubsection{}\label{It:3grow} $\sigma = (0, 1, \dots, s_1-1, s_1 + 1, 
\dots, s_2 - 1, s_2 + 1, \dots, s_3 - 1, s_3 + 1,\dots, n-1, s_1, s_2, 
s_3)$ where $s_1 < s_2 < s_3$ and $s_3 \le n-4$. In this case $u_i(\sigma) 
= 0$ for $i \le n-4$, $u_{n-3}(\sigma) = s_1 + 1$, $u_{n-2}(\sigma) = s_2$, 
and $u_{n-1}(\sigma) = s_3 - 1$. Consider a pair $(k,l)$ where $l_i = 0$ 
for $i \le n-4$, $l_{n-3} = s_1 + 1$, $l_{n-2} = s_2$, and $l_{n-1} = s_3 - 
1$. It is easy to see that $(k,l)$ is admissible  and $f(k,l) = 
p(\sigma,k)$. Vice versa, if $(k,l)$ is an admissible pair with $l$ such 
that $l_0 = \dots = l_{n-4} = 0$ and $1 \le l_{n-3} \le l_{n-2} \le l_{n-1} 
\le n-2$, one can construct an HL-pair $(\sigma,k)$ by the formula above 
with $s_1 = l_{n-3}-1$, $s_2 = l_{n-2}$, and $s_3 = l_{n-1} + 1$. 

\subsubsection{}\label{It:2inv} $\sigma = (0, 1, \dots, s_2-1, s_2 + 1, 
\dots, s_1 - 1, s_1 + 1, \dots, n-1, s_1, s_2)$ where $s_2 < s_1 \le n-2$. 
In this case $u_i(\sigma) = 0$ for $i \le n-3$, $u_{n-2}(\sigma) = s_1$ and 
$u_{n-1}(\sigma) = n-1$. Consider a pair $(k,l)$ where $l_i = 0$ for $i \le 
n-3$, $l_{n-2} = s_1$ and $l_{n-1} = s_2$. It is easy to see that $(k,l)$ 
is admissible and $f(k,l) = p(\sigma,k)$. Vice versa, if $(k,l)$ is an 
admissible pair with $l$ such that $l_0 = \dots = l_{n-3} = 0$ and $0 \le 
l_{n-2} > l_{n-1} \ge 0$, one can construct an HL-pair $(\sigma,k)$ by the 
formula above with $s_1 = l_{n-2}$ and $s_2 = l_{n-1}$. 

So, it follows from the conjecture of \cite{HL} that the dimension of 
$H_{a,3}$ equals the number of admissible pairs with $\sum_i k_i = a$ and 
$l$ either like in Case \ref{It:3grow} or like in Case \ref{It:2inv}. This 
coincides with the conjecture of \cite{BS} describing the set $b^{-1}(3)$.

\subsection{} Let $\name{maj}(\sigma) = 4$. For this case \cite{BS} 
does not provide a conjecture, and the following comparison with \cite{HL} 
explains the nature of the difficulties. The permutation $\sigma$ now has 
one of the three forms:

\subsubsection{} $\sigma = (0, 1, \dots, s_1-1, s_1 + 1, \dots, s_2 - 1, 
s_2 + 1, \dots, s_3 - 1, s_3 + 1,\dots, s_4 - 1, s_4 + 1,\dots, n-1, s_1, 
s_2, s_3)$ where $s_1 < s_2 < s_3 < s_4$ and $s_4 \le n-5$. This case is 
analogous to the case \ref{It:3grow}. HL-pair $(\sigma,k)$ with such 
$\sigma$ can be related to an admissible pair $(k,l)$ where $l_i = 0$ for 
$i \le n-5$ and $l_{n-4} = s_1 + 1$, $l_{n-3} = s_2$, $l_{n-2} = s_3 - 1$, 
$l_{n-1} = s_4 - 2$; then $f(k,l) = p(\sigma,k)$.

\subsubsection{} $\sigma = (0, 1, \dots, s_1-1, s_1 + 1, \dots, s_3 - 1, 
s_3 + 1, \dots, s_2 - 1, s_2 + 1,\dots, n-1, s_1, s_2, s_3)$ where $s_1 < 
s_3 < s_2$ and $s_1 \le n-3$. In this case $u_i(\sigma) = 0$ for $i \le 
n-4$, $u_{n-3}(\sigma) = s_1 + 1$, $u_{n-2}(\sigma) = s_2$ and 
$u_{n-1}(\sigma) = n-1$. An HL-pair $(\sigma,k)$ can be related to an 
admissible pair $(k,l)$ where $l_i = 0$ for $i \le n-4$, $l_{n-3} = s_1 + 
1$, $l_{n-2} = s_2$ and $l_{n-1} = s_3$; then $f(k,l) = p(\sigma,k)$.

\subsubsection{} $\sigma = (0, 1, \dots, s_3-1, s_3 + 1, \dots, s_1 - 1, 
s_1 + 1, \dots, s_2 - 1, s_2 + 1,\dots, n-1, s_1, s_2, s_3)$ where $s_3 < 
s_1 < s_2$ and $s_1 \le n-3$. This is the case where difficulties arise. 
Here $u_i(\sigma) = 0$ for $i \le n-4$, $u_{n-3}(\sigma) = s_1 + 1$, 
$u_{n-2}(\sigma) = s_2$ but, unlike the previous case, $u_{n-1}(\sigma) = 
n-2$. There is no correspondence $\sigma \mapsto l$ such that if 
$(\sigma,k)$ is an HL-pair of this form then $(k,l)$ is an admissible pair 
with $f(k,l) = p(\sigma,k)$; nor I was able to find a more refined 
correspondence which changes $k$ accordingly.

\end{document}